\documentclass[10pt]{elsart3-1} 
\usepackage{amssymb}
\usepackage[english,francais]{babel} 
\usepackage[figuresright]{rotating}
\usepackage{amsmath, amsfonts}

\usepackage[most]{tcolorbox} 

\oddsidemargin -.in
\newtheorem{e-proposition}[theorem]{Proposition} 

\newtheorem{e-definition}[theorem]{Definition\rm}

\newtheorem{theoreme}{Th\'eor\`eme}[section]

\newtheorem{proposition}[theoreme]{Proposition}

\newtheorem{definition}[theoreme]{Definition\rm}

\renewcommand{\theequation}{\arabic{equation}}
\numberwithin{equation}{section}
\numberwithin{figure}{section}

\newcommand \bse {\begin{subequations}}
\newcommand \ese {\end{subequations}}

\def\og{\leavevmode\raise.3ex\hbox{$\scriptscriptstyle\langle\!\langle$~}}
\def\fg{\leavevmode\raise.3ex\hbox{~$\!\scriptscriptstyle\,\rangle\!\rangle$}}

\renewcommand \sharp {{}}
\newcommand \Deltax h
\newcommand \la \langle
\newcommand \ra \rangle

\newcommand \Dcal {\mathcal D}
\newcommand \RN {{{\mathbb R}^N}}
\newcommand \Rd {{\mathbb R}^d}

\newcommand \Ucal {\mathcal U} 
\newcommand \Ucals {\mathcal U_\sharp} 
\newcommand \us {{u_\sharp}} 

\newcommand \Hs {{H_\sharp}} 
\newcommand \Ds {{D_\sharp}} 
\newcommand \vs {v_\sharp}
\newcommand \Ss {S_\sharp}

\newcommand \fs {{f_\sharp}{}}
\newcommand \Bs {{B_\sharp}{}}

\newcommand \Us {U_\sharp}
\newcommand \Ks {K_\sharp}
\newcommand \Ls {L_\sharp}
\newcommand \Fs {F_\sharp{}}

\newcommand \hs {h_\sharp}

\newcommand \bs {b_\sharp}
\newcommand \fsTot {f^{\sharp \text{Tot}}}

\newcommand \fsDiff {f^{\sharp \text{Diff}}}
\newcommand \fsDisp {f^{\sharp \text{Disp}}}
\newcommand \FsTot {F^{\sharp \text{Tot}}}
\newcommand \FsDiff {F^{\sharp \text{Diff}}}
\newcommand \FsDisp {F^{\sharp \text{Disp}}}
\newcommand \DsDiff {D^{\sharp \text{Diff}}}
\newcommand \bei {\begin{itemize}}
\newcommand \eei {\end{itemize}}  
\newcommand \RR   {\mathbb R}
\newcommand \be  {\begin{equation}}
\newcommand \ee  {\end{equation}}
\newcommand \del \partial
\newcommand \eps \varepsilon
  
\newcommand \bel {\be\label}

\journal{the Acad\'emie des sciences}
\begin{document} 
\centerline{}
\begin{frontmatter}
 
\vskip-2.5cm

\selectlanguage{english}
\title{Augmented hyperbolic models with diffusive-dispersive shocks}

\selectlanguage{english}
\author[authorlabel1]{Philippe G. LeFloch}
\ead{contact@philippelefloch.org}
\and 
\author[authorlabel2]{Allen M. Tesdall}
\ead{allen.tesdall@csi.cuny.edu}

\address[authorlabel1]{Laboratoire Jacques-Louis Lions \& Centre National de la Recherche Scientifique 
\\
Sorbonne Universit\'e, 4 Place Jussieu, 75252 Paris, France. } 

\address[authorlabel2]{ Department of Mathematics, City University of New York, College of Staten Island, and Physics Program, The Graduate Center, City University of New York, New York, U.S.A. 
\\
November 2019}

\medskip

\begin{abstract}
\selectlanguage{english}   
Given a first-order nonlinear hyperbolic system of conservation laws  
endowed with a convex entropy-entropy flux pair, we can consider the class of weak solutions containing shock waves depending upon some small scale parameters. In this Note, we define and derive several classes of entropy-dissipating augmented models, as we call them, which involve (possibly nonlinear) second- and third-order augmentation terms. Such terms typically arise in continuum physics and model the viscosity and capillarity effects in a fluid, for instance. By introducing a new notion of positive entropy production that concerns general functions (rather than solutions) we can easily check the entropy-dissipating property for a broad class of augmented models. The weak solutions associated with the zero diffusion/dispersion limit may contain (nonclassical undercompressive) shocks whose selection is determined from these diffusive and dispersive effects (for instance by using traveling wave solutions), and having a classification of the models, as we propose, is essential for developing a general theory. 
%
%
%
%
%
%
%
%
%

\end{abstract}
\end{frontmatter}

%
   

\vskip-.65cm

\selectlanguage{english} 



\section{Small-scale sensitive shocks and positive entropy production}
 
\paragraph*{Continuum physics modeling from small-scales to macro-scales} 

Many models in continuum physics involve augmentation terms containing small parameters such as the viscosity, capillarity, heat conduction, or other parameters related to relaxation effects, the Hall effect, etc.  
These terms are typically modeled by second- or third-order derivatives which are taken into account in the fundamental conservation principles of continuum physics. For instance in a compressible fluid, when the viscosity effects are dominant, no region with a sharp gradient can form. On the other hand, when the capillarity effects are dominant, highly oscillating patterns are observed near sharp gradients of the solutions. In the present Note, the regime we are interested in is the one when the viscosity and capillarity effects are small while being kept {\sl in balance} with each other \cite{LeFloch-Freiburg}. It is natural to assume that the 
orders of magnitude of the physical coefficients (which depend on the fluid or solid material under consideration) arising in (for instance) a combination of viscosity and capillarity terms like (see below)
$
\nu \, u_{xx} + \kappa \, u_{xxx} 
$
are such that the ratio 
$
\xi ={\kappa / \nu^2}$ 
is of order $1$. We can then write these terms (for some $\eta, \eps$) as 
$
\nu \, u_{xx} + \kappa \, u_{xxx} 
=
 \eta \eps u_x^\eps + \xi \eta^2 \eps^2 u_{xx}^\eps. 
$

We are still assuming that the hyperbolic aspects of the flow are dominant, and the visocity and capillarity effects are relevant only in a neighborhood of large gradients of the solutions. In the theory of nonclassical shocks one is interested in describing the macro-scale features, that is, only the {\sl limit $\eps \to 0$} 
and to do so we need to extract some ``information''  from the small-scales. Importantly, for the present Note, the global dynamics of the shocks turn out to depend upon the small-scale physical modeling, and it is our aim here to provide a framework in which general classes of interest can be derived and analyzed.


\paragraph*{Systems of conservation laws endowed with an entropy}

We are interested in systems of the form 
\bel{eq:101}
\del_t u + \sum_{j=1}^d \del_j f^j(u) = 0, 
\qquad u=u(t,x) \in \Ucal \subset \RN, \quad t \geq 0, \, x \in \Rd, 
\ee
with unknown $u=(u_a)=(u_1, \ldots, u_N) \in \Ucal$ (an open subset of $\RN$ containing $0$), 
and we assume such a system to be endowed with a  convex
 entropy-entropy flux pair $(U, F)$ satisfying, by definition, 
$D^2 U > 0$ and $DF^j = (D U)^T D f^j$, 
and normalized so that $U(0) = 0$, as well as $f(0)=0$ and $F(0)=0$. 
By definition, physically meaningful solutions, also called entropy solutions, must satisfy the entropy inequality 
\bel{eq:102}
\del_t U(u) +  \sum_{j=1}^d \del_j F^j(u)  \leq 0, 
\qquad 
t \geq 0, \, x \in \Rd.
\ee

The so-called entropy variable defined by 
\be
\aligned
& v = \vs(u) := \nabla U(u) \in \Ucals := \vs(\Ucal), 
\qquad 
\quad
 u = \us(v) := \big(\nabla U\big)^{-1}(v), 
\endaligned
\ee
will play a fundamental role in our approach. Performing the change of variable $u \in \Ucal \mapsto v\in \nabla U(\Ucal)\subset \RN$ and setting 
\be
\fs^j(v) := f^j(u), 
\quad
\Us(v) := U(u), \quad 
\Fs^j(v) := F^j(u),
\ee
we rewrite \eqref{eq:101} and \eqref{eq:102} in the form: 
\bel{eq:103}
\del_t \us(v) +  \sum_{j=1}^d \del_j \fs^j(v) = 0, 
\qquad  
\del_t \Us(v) +   \sum_{j=1}^d \del_j \Fs^j(v)  \leq 0,
\qquad 
t \geq 0, \, x \in \Rd.
\ee 
While many systems in fluid and solid dynamics fit into this class, still it does not fully describe the dynamics of small-scale sensitive shocks.


\paragraph*{A broad class of augmented systems with diffusion and dispersion}

More relevant physical models read 
\bel{eq:101-general}
\del_t u^\eps + \sum_{j=1}^d \del_j f^j(u^\eps) 
=  \sum_{j=1}^d \del_j S^j[u^\eps],
\qquad u^\eps=u^\eps(t,x),  
\ee
in which $\eps>0$ is a small parameter representing the amount of viscosity, capillarity, heat conduction, etc. in the fluid or material under consideration
and $S^j[u^\eps]$ depends on (suitably scaled)
first- and higher-order derivatives $\eps \del_j u^\eps$ and $\eps^2 \del_j \del_k u^\eps$. 
We assume the normalization that if $u^\eps=u^\eps(x)$ is a constant function then  $S^j[u^\eps]$ vanishes identically. 

For a physically realistic model, we expect certain positivity and structural conditions to be imposed, especially the entropy inequality \eqref{eq:102}. 
As we will see, some terms may add (non-negative) contributions
to the entropy itself, while some terms will add (non-negative) contributions to the entropy dissipation, while other terms with divergence form will not contribute to the total entropy balance. 

The fundamental requirement we want to impose on the higher-order terms is the following {\bf sign condition}:
\bel{eq:sign}
\limsup_{\eps \to 0}  \iint_{\RR^+\times \RR^N} 
\sum_{j=1}^d 
\nabla U(u^\eps)^T  \del_j S^j[u^\eps] \, \theta dx dt
\leq 0
\ee
for {\sl every solution} $u^\eps$ to \eqref{eq:101-general} and every (smooth and compactly supported) test-function $\theta=\theta(t,x) \geq 0$. 
Clearly, if the condition \eqref{eq:sign} holds, then we deduce from 
the augmented model \eqref{eq:101-general} that the limiting solution $u:=\lim u^\eps$ (if it exists in a suitable functional space) satisfies the entropy inequality 
\eqref{eq:102}. However, the condition \eqref{eq:sign} is {\sl not explicit} enough to be useful in practice, for instance for numerical discretization. We are going to present suitable classes of models that are, both, physically relevant and numerically tractable, and enjoy a  {\bf positive entropy production property}. Our notion guides us in identifying the interesting classes of models (and later designing the schemes adapted to these models). 

It is important to handle the condition \eqref{eq:sign} in the entropy variable, that is, to write
\bel{eq:103-with-v}
\del_t \us(v^\eps) +  \sum_{j=1}^d \del_j \fs^j(v^\eps) 
= \sum_{j=1}^d \del_j \Ss^j[v^\eps],
\qquad v^\eps=v^\eps(t,x), 
\qquad 
t \geq 0, \, x \in \Rd, 
\ee
together with the sign condition 
\bel{eq:sign-entro}
\limsup_{\eps \to 0} \iint_{\RR^+\times \RR^N} 
\sum_{j=1}^d 
(v^\eps)^T  \del_j \Ss^j[v^\eps] \, \theta dx dt
\leq 0. 
\ee 
Here, we have set $\Ss^j[v^\eps] :=  S^j[u^\eps]$
and \eqref{eq:sign-entro} is required for all test-functions $\theta=\theta(t,x) \geq 0$ and all solutions $v^\eps=v^\eps(x)$ to \eqref{eq:103-with-v}.


\section{The positive entropy production property for augmentation terms}
\label{sec:generalpositivity} 

\subsection{Arbitrary functions as test-functions}

We proceed by suppressing the time integral in the dissipation bound that arises from \eqref{eq:sign} 
and we propose the following notion which, importantly, {\sl no longer refers} to the PDE under consideration, but imposes a condition on the augmentation terms for {\sl general functions} rather than solutions. 

\begin{definition}
\label{def=main}
A nonlinear expression $S[w^\eps] = (S^j[w^\eps])_{1 \leq j \leq d}$ of a sequence of functions $w^\eps: \Rd \to \Ucals \subset \RN$ 
involving $w^\eps$ and its (rescaled) derivatives $\eps \del_j w^\eps$ and $\eps^2 \del_j \del_k w^\eps$ 
has the {\bf positive entropy production property} if
\bel{eq:sign-entro-bad}
\liminf_{\eps \to 0} 
\int_{\Rd}   
\nabla w^\eps \cdot  \Ss[w^\eps] \, \theta dx \geq 0
\ee
(with $\nabla w^\eps \cdot \Ss[w^\eps] :=  \sum_{j=1}^N \del_j w_\eps^T \Ss^j[w^\eps]$) 
for all test-functions $\theta=\theta(x) \geq 0$
and for any sequence of bounded functions $w^\eps=w^\eps(x)$ with bounded 
total dissipation
\bel{eq:total-dissip-3} 
\limsup_{\eps \to 0} \| w^\eps\|_{L^\infty(\RN)} 
+
\limsup_{\eps \to 0} \int_{\Rd} \eps \,  |\nabla w^\eps|^2 dx 
< +\infty. 
\ee
\end{definition}

We can prove that if the augmentation terms satisfy the positive entropy production property, then 
the dissipation bound in \eqref{eq:total-dissip-3} follows from the sole assumption that the total entropy is bounded. More generally, the augmentation terms could also depend on time-derivatives, but such a generalization is more involved since the time variable must be handled by using the equation \cite{PLF-AMT-1}. 

Consider any nonlinear expression $S[w^\eps] = (S^a[w^\eps])_{1 \leq a \leq d}$,  
having the positive entropy production property (cf.~Definition~\ref{def=main})
and depending upon a sequence of functions $w^\eps: \Rd \to \RN$ and its (rescaled) derivatives $\eps \del_a w^\eps$
and $\eps^2 \del_a \del_b w^\eps$ with $a,b= 1, \ldots, d$.  
Assuming the usual bound 
\bel{eq:total-dissip-3-bad} 
\limsup_{\eps \to 0} \| w^\eps\|_{L^\infty(\RN)} 
+
\limsup_{\eps \to 0} \int_{\Rd} \eps \,  |\nabla w^\eps|^2 \, dx 
< +\infty, 
\ee
we can associate to $w=(w^\eps)$ a locally bounded measure $\mu_w$ defined over $\RN$ so that 
for every test-function $\theta=\theta(x) \geq 0$ 
\bel{eq:sign-entro-2-bad}  
\la \mu_w, \theta \ra 
:= 
\liminf_{\eps \to 0} 
\int_{\Rd}  \eps \, 
\nabla w^\eps \cdot  \Ss[w^\eps] \, \theta dx \geq 0. 
\ee  
This measure depends upon the choice of the sequence. As long as classical shocks are concerned, the positivity of this measure implies that the physically meaningful shock wave is selected by our augmented model. 
However, the actual values of this measure are required for selecting of nonclassical (undercompressive) shocks.


\subsection{Augmented models with space derivatives}
\label{sec:onlyspace}

While our results hold in several space dimensions, they are easier to present in one space dimension, so we restrict attention here to 
\bel{eq:101-one}
u_t + f(u)_x = 0, \qquad u=u(t,x) \in \Ucal \subset \RN, 
\ee
with unknown $u=(u_a)=(u_1, \ldots, u_N) \in \Ucal$ (an open subset of $\RN$). 
We assume this system to be endowed with a  convex  entropy-entropy flux pair $(U, F)$, as described in the introduction, 
and we recall the entropy inequality 
\bel{eq:102-one}
U(u)_t + F(u)_x \leq 0. 
\ee
In the entropy variable $v= \vs(u) := \nabla U(u)$ 
after performing the change of variable $u \in \Ucal \mapsto v\in \nabla U(\Ucal)\subset \RN$ 
we have  
\bel{eq:103-deux}
\us(v)_t + \fs(v)_x = 0, 
\ee
and observing that 
\bel{eq:8418}
\aligned
v^T \us(v)_t  = \Us(v)_t, 
\qquad 
v^T \fs(v)_x  & = \Fs(v)_x, 
\endaligned
\ee
we obtain 
\bel{eq:102-VG}
\Us(v)_t + \Fs(v)_x \leq 0. 
\ee 
As already mentioned, it is natural to express the high-order contributions in terms of $v$ and we now proceed 
by listing classes of augmented models of increasing difficulty. 
For simplicity in the presentation, we often suppress the subscript $\eps$ except when emphasis is necessary. 

It will be useful to have a short-hand notation for the ``total'' flux and entropy flux, that is, we will write 
the augmented model in the form
\bse
\be
\us[v]_t + \fsTot[v]_x  = 0, 
\ee
and its entropy balance law in the form 
\be
\Us[v]_t + \FsTot[v]_x = \Ds[v], 
\ee
in which the following expressions can be determined for each model of interest: 
\be
\fsTot[v] = \fs(v) + \fsDiff[v] + \fsDisp[v], 
\qquad 
\FsTot[v] = \Fs(v) + \FsDiff[v] + \FsDisp[v]. 
\ee
\ese


\subsection{Linear diffusion and linear dispersion} 

We begin with augmentation terms that are linear in $v_x$ and $v_{xx}$ with constant coefficients. 
Namely, we consider
\bel{eq:109-constant} 
\aligned
\us(v)_t + \fs(v)_x 
 = \eps \, \Bs v_{xx} + \eps^2 \Ks v_{xxx} 
= \big( \Ss[v^\eps] \big)_x, 
\endaligned
\ee
in which $\Bs$ and $\Ks$ are given $(N\times N)$-matrices. Our sign condition \eqref{eq:sign-entro-bad}
is equivalent to saying 
\bel{eq:266}
\liminf_{\eps \to 0} \int_{\RR}
(w^\eps_x)^T \Big(
 \eps \, \Bs w^\eps_{x} + \eps^2 \Ks w^\eps_{xx} 
\Big)
\, \theta dx 
\geq 0
\ee
for all test-functions $\theta=\theta(x) \geq 0$
and all sequences  $w^\eps: \Rd \to \Ucals \subset \RN$ satisfying \eqref{eq:total-dissip-3}. 

\begin{proposition}[Linear diffusion and linear dispersion]
\label{propo-347}
A \underline{necessary condition} for the one-dimensional augmented model \eqref{eq:109-constant} to satisfy the positive entropy production property in Definition~\ref{def=main}
(see \eqref{eq:266}) is: 
\be
\aligned
&  \Bs^T + \Bs  \text{ is a non-negative matrix,} 
\qquad
 \Ks    \text{ is a symmetric matrix.} 
\endaligned
\ee
Moreover, the entropy inequality holds as follows: 
\bse
\bel{eq:109-constant-10} 
\aligned
& \Us(v)_t + \FsTot[v]_x
 = \Ds[v] = -  {1 \over 2}\eps v_x^T \big( \Bs + \Bs^T \big) v_x \leq 0,  
\endaligned
\ee
where 
\bel{eq:109-constant-10b} 
\aligned
& \FsTot[v] = \Fs(v) + \FsDiff[v] + \FsDisp[v], 
\hskip2.cm
 \FsDiff[v] = -  \eps v^T \Bs v_x, 
\\
& \FsDisp[v] := {\eps^2 \over 2} 
\Big( 3 v_x^T \Ks v_x  -  \big( v^T \Ks v \big)_{xx}
\Big). 
\endaligned
\ee 
\ese 
\end{proposition}


\subsection{Nonlinear diffusion and nonlinear dispersion. Sufficient conditions} 

We proceed by successive generalizations of the augmentation terms and we now consider  a nonlinear version of \eqref{eq:109-constant}, that is, 
\bel{eq:109-2} 
\aligned
& \us(v)_t + \fs(v)_x 
 = \eps \, \big( \Bs(v) \,v_x \big)_x  
+ \eps^2 \, \Big( \nabla\hs(v)^T H \hs(v)_{xx} \Big)_x = \eps \big( \Ss[v^\eps] \big)_x, 
\endaligned
\ee
in which $\Bs=\Bs(v)$ is a given $(N\times N)$-matrix-valued mapping 
and $\hs=\hs(v)$ is an $N$-vector-valued mapping, while $H$ is a constant matrix. The choice of this
 structure will be further motivated below; at this stage, we consider \eqref{eq:109-2} 
as an interesting broad class of models. 

Our sign condition \eqref{eq:sign-entro-bad} is equivalent to saying 
\bel{eq:266-variable}
\limsup_{\eps \to 0} \int_{\RR}
(w^\eps_x)^T \Big(
 \eps \, \Bs (w^\eps) w^\eps_x + \eps^2 \nabla\hs(w^\eps)^T H \hs(w^\eps)_{xx}
\Big)
\, \theta dx 
\geq 0
\ee
for all test-functions $\theta=\theta(x) \geq 0$ and all sequences $w^\eps: \Rd \to \Ucals \subset \RN$
 satisfying \eqref{eq:total-dissip-3}. 
It is convenient to introduce the notation 
$\us(v)_t + \big( \fsTot[v] \big)_x = 0$ 
where we now have  
\bse
\be
\aligned
\fsTot[v] 
& = 
\fsTot(v, \eps v_x, \eps^2 v_{xx}) 
 := \fs(v) + \fsDiff(v, \eps v_x)  + \fsDisp(v, \eps v_x, \eps^2 v_{xx}), 
\endaligned
\ee 
with 
\be
\fsDiff[v]
= 
\fsDiff(v, \eps v_x) := \eps \Bs(v) \, v_x, 
\ee
and 
\be
\aligned
\fsDisp[v) 
& = 
\fsDisp(v, \eps v_x, \eps^2 v_{xx})  
:= -\eps^2 \nabla\hs(v)^T H \, \hs(v)_{xx}   
\\
& =  - \eps^2 \nabla\hs(v)^T H \,  \nabla\hs(v) v_{xx} - \eps^2 \nabla\hs(v)^T H \,  \nabla^2\hs(v) (v_x, v_x). 
\endaligned
\ee 
\ese 

\begin{proposition}[Nonlinear diffusion and nonlinear dispersion]
\label{Prop:345} 
Provided the following properties hold: 
\bel{eq:294}
\aligned
&  \Bs(v)^T + \Bs(v)  && \text{ is a non-negative matrix,} 
\\
& H             && \text{ is a symmetric matrix,} 
\endaligned
\ee
the augmented model \eqref{eq:109-2}  satisfies the positive entropy production property  in Definition~\ref{def=main} 
(see \eqref{eq:266-variable}) 
and, furthermore, the entropy inequality holds as follows:  
\bse
\bel{eq:109-constant-10-nl} 
\aligned
& \Us(v)_t + \big( \FsTot[v] \big)_x
= \Ds[v] 
= -  {1 \over 2}\eps v_x^T \big( \Bs (v)+ \Bs^T(v) \big) v_x \leq 0, 
\endaligned
\ee
where 
\be
\aligned 
&
\FsTot[v] 
= \Fs(v) + \FsDiff[v] + \FsDisp[v], 
\qquad
\hskip1.cm 
\FsDiff[v] 
=  - \eps  v^T \Bs(v) v_x, 
\\
& 
\FsDisp[v] 
:=   -  v^T\nabla\hs(v)^T \, \Hs \hs(v)_{xx} + {1\over 2} \hs(v)_x H \hs(v)_x  
\\
& 
=  {\eps^2 \over 2} 
\Big( 3 v_x^T \Ks(v) v_x  -  \Big( v^T \Ks (v)v \Big)_{xx} \Big) + \eps^2 \nabla^2 \hs (v, v_x)^T \Hs \hs(v)_x 
=  -  {\eps^2 \over 2} \big( v^T \Ks (v)v \big)_{xx} + \eps^2 v_x^T  \Ls(v) v_x, 
\endaligned
\ee 
and 
\be
\aligned
\Ks(v) & := \nabla\hs(v)^T \Hs \nabla\hs(v),
\qquad
\Ls(v) := {3 \over 2} \Ks(v) + \nabla^2 \hs (v,  \cdot)^T H \nabla \hs(v).
\endaligned
\ee
\ese 
\end{proposition}

Observe that $\Ks(v)$ is automatically a symmetric matrix when $H$ is symmetric. 
Note also that the model \eqref{eq:109-2} reduces to \eqref{eq:109-constant} if $\hs(v)$ is chosen to be $v$. 
Observe that $\FsDisp[v]$ is the sum of a term in a divergence form and a quadratic term, while  
$\FsDiff[v] = - \eps  v^T \Bs(v) v_x$ would also have a divergence form if we restrict attention to diffusion matrices deriving from a 
scalar potential 
$\bs$ of the form $v^T \Bs = \nabla \bs$ and thus $\FsDiff[v] = - \eps \bs(v)_x$. 
In the applications, the quadratic term related to $\Bs(v)^T + \Bs(v)$ will indeed have the required sign, while 
the dispersion contribution from $\Ls(v)$ will also often be non-negative. 

 
\section{Nonlinear diffusion and nonlinear dispersion}

\subsection{Necessary conditions}

We will now discuss the question of whether the term $\nabla\hs(v)^T H \hs(v)_{xx}$ is the most general diffusive term ensuring a good entropy structure. 
For clarity in the discussion we restrict attention to scalar equations with a quadratic entropy, that is, with $F' := u \, f'$ 
\bel{eq:gf99}
\aligned
u_t + f(u)_x & = 0, 
\qquad
{1 \over 2} (u^2)_t + F(u)_x  \leq 0. 
\endaligned
\ee

\begin{proposition}[ Necessary conditions. I]
 Consider the conservation law with quadratic entropy \eqref{eq:gf99}. 
Given some continuous functions $B = B(u)$, $k_1=k_1(u)$, and $k_2=k_2(u)$, the augmentation terms in the diffusive-dispersive model 
\be
u_t + f(u)_x = 
\eps (B(u) u_x)_x + \eps^2 \big( k_1(u) (k_2(u) u_x)_x \big)_x 
= \big( S[u, \eps u_x, \eps^2 u_{xx}]\big)_x 
\ee
satisfy the positive entropy dissipation property in Definition~\ref{def=main} 
if and only if the functions satisfy 
\be
B(u) \geq 0, \qquad
k_1(u) = c \, k_2(u), \qquad u \in \RR, 
\ee
where $c$ is an arbitrary constant. This is precisely the structure already analyzed in Proposition~\ref{Prop:345}. 
\end{proposition}


We complete our discussion by analyzing the most general diffusive term that is linear with respect to the highest derivative. An example ensuring the favorable signs below is 
$S_1(u,v) = b_1(u) \, v^{2p+1} $
and 
$S_2(u,v) = b_2(u) \, v^{2q+1}$ provided $b_1 \geq 0$ and $b_2' \leq 0$. 

\begin{proposition}[ Necessary conditions. II]
 Consider the conservation law with quadratic entropy \eqref{eq:gf99}. 
Given some functions $S_1 = S_1 (u, \eps u_x)$ and $S_2 = S_2(u, \eps u_x)$, assumed to be analytic in their arguments, 
the augmentation terms in the diffusive-dispersive model 
\bel{eq:yet-a-newmodel}
u_t + f(u)_x = \big( S_1 (u, \eps u_x) + \eps^2 u_{xx} S_2 (u, \eps u_x) \big)_x
\ee
satisfy the positive entropy dissipation property (in Definition~\ref{def=main})
if and only if  
\bel{eq:HD773} 
\aligned 
v S_1 (u, v) - v^3 \, \del_1 \overline  S_2 \big(u, v) \geq 0, \qquad u, v \in \RR, 
\endaligned
\ee  
where 
\be
\overline  S_2(u,v)
 := \int_0^1 \int_0^1 S_2 \big(u, v (1 + m (s-1)) \big) \,  (1-s) \, ds dm.
\ee
The corresponding entropy balance law then reads 
\bse
\be 
\aligned
{1 \over 2} (u^2)_t + \big( \FsTot[u] \big)_x
& = \Ds[u] 
  = - u_x S_1 (u, \eps u_x) + \eps^2 u_x^3 \del_1 \overline  S_2 \big(u, \eps u_x)
 \leq 0, 
\endaligned
\ee
where 
\be
\aligned  
&  \FsDiff[v] := u S_1 (u, \eps u_x), 
\qquad
 \FsDisp[v] := - \eps^2 u_{xx} u S_2 (u, \eps u_x) - \eps^2 u_x^2 \overline S_2 (u, \eps u_x). 
\endaligned
\ee 
\ese
\end{proposition} 


\subsection{Global entropy balance law}

Denoting now the solution by $v^\eps$ and integrating in space over $\RR$, we find (provided the solution decays appropriately at infinity)
$$
{d \over dt} \int_\RR \Us(v^\eps)  \, dx 
= \int_\RR \DsDiff(v^\eps, \eps v^\eps_x) \, dx 
\leq 0,
$$
and in particular 
$
t \geq 0 \mapsto \int_\RR \Us(v^\eps(t,x)) \, dx$
is non-increasing. 
We conclude our discussion in the present section with the following observation. 
Note that the diffusion and dispersion coefficients and matrices may well be degenerate.  Finally, we make the following important observation. For every model we introduced, 
if the sequence of initial data $v^\eps(0, \cdot)$ has uniformly bounded entropy and the solutions $v^\eps$ 
remain uniformly bounded in sup-norm, i.e. 
$
\limsup_{\eps \to 0} \| v^\eps\|_{L^\infty(\RR)} < +\infty
$
and converge almost everywhere to some limit $v$, i.e.
$
v := \lim_{\eps \to 0} v^\eps, 
$
then 
one has 
$$
\aligned
& \fsTot[v] - \fs(v) \to 0 \, \text{ and } \, 
\FsTot[v] - \Fs(v) \to 0 \quad \text{ in the sense of distributions,}
\\
& \Dcal[v^\eps] 
\text{ is a sequence of locally uniformly bounded non-negative measures,}
\endaligned
$$
and, moreover, the limit is a weak solution satisfying the entropy inequality, that is, 
$$
\aligned
\us(v)_t + \fs(v)_x & = 0,
\qquad
\Us(v)_t + \Fs(v)_x  \leq 0.
\endaligned
$$
We  also emphasize that specific models from nonlinear elasticity and phase transition dynamics are found to fit within our setting. 

For the proofs of the results announced in this Note, we refer to \cite{PLF-AMT-1}. For the equations derived here, we have also developed adapted numerical schemes that are ``structure-preserving'' \cite{PLF-AMT-2}. For nonlinear hyperbolic problems, many strategies were proposed in recent years in order to discretize certain algebraic or differential properties; among many contributions, see \cite{BoscarinoRS} and \cite{Chertock}. A distinct feature of diffusive-dispersive shocks is the richer variety of waves that are observed; yet these standard structure-preserving techniques are relevant for dealing with the more involved models introduced in the present Note. 


\small 

\paragraph*{\bf Acknowledgments.}
 
This paper was completed while the first author (PLF) was a visiting fellow at the Courant Institute of Mathematical Sciences at New York University during the Academic year 2018--2019. 
The work of the second author (AMT) was supported by NSF under grant DMS-1516131, and by CUNY Research Foundation under PSC-CUNY Grant \#60365-00 48. 
The first author (PLF) also is grateful to AMT for his financial support through his NSF and CUNY grants. 


\bibliographystyle{plain}

\end{document}